\documentclass[letterpaper, 10 pt, conference]{ieeeconf}  

\IEEEoverridecommandlockouts                              

\usepackage{graphicx} 
\usepackage{amsmath} 
\usepackage{amsbsy}
\usepackage{amssymb}  
\usepackage{cleveref}
\usepackage{mathrsfs}
\usepackage{float}
\usepackage[bottom]{footmisc}

\DeclareMathOperator*{\argmin}{argmin}

\newcommand{\R}{\mathbb R}
\newcommand{\bx}{\boldsymbol{x}}
\newcommand{\by}{\boldsymbol{y}}
\newcommand{\bs}{\boldsymbol{\theta}}
\newcommand{\bv}{\boldsymbol{v}}
\newcommand{\bw}{\boldsymbol{\omega}}
\newcommand{\w}{\omega}
\newcommand{\abs}[1]{\left \lvert #1 \right \rvert}

\graphicspath{{images/}}

\title{\LARGE \bf
Time-Optimal Paths for Simple Cars with Moving Obstacles in the Hamilton-Jacobi Formulation
}

\author{Christian Parkinson$^{1}$ and Madeline Ceccia$^{2}$
\thanks{$^{1}$Christian Parkinson is a postdoctoral research associate with the Department of Mathematics, 
        University of Arizona, 617 N. Santa Rita Ave, Tucson, AZ, 85721
        {\tt\small chparkin@math.arizona.edu}}%
\thanks{$^{2}$Madeline Ceccia is a student in the Department of Mathematics, California State University - Fullerton, 800 North State College Blvd.
Fullerton, CA 92831
        {\tt\small madelinececcia@csu.fullerton.edu}}%
}

\begin{document}

\maketitle
\thispagestyle{empty}
\pagestyle{empty}

\begin{abstract}

We consider the problem of time-optimal path planning for simple nonholonomic vehicles. In previous similar work, the vehicle has been simplified to a point mass and the obstacles have been stationary. Our formulation accounts for a rectangular vehicle, and involves the dynamic programming principle and a time-dependent Hamilton-Jacobi-Bellman (HJB) formulation which allows for moving obstacles. To our knowledge, this is the first HJB formulation of the problem which allows for moving obstacles. We design an upwind finite difference scheme to approximate the equation and demonstrate the efficacy of our model with a few synthetic examples. 

\end{abstract}

\section{INTRODUCTION}

\begin{figure}[b!]
\centering
\includegraphics[width=0.35\textwidth]{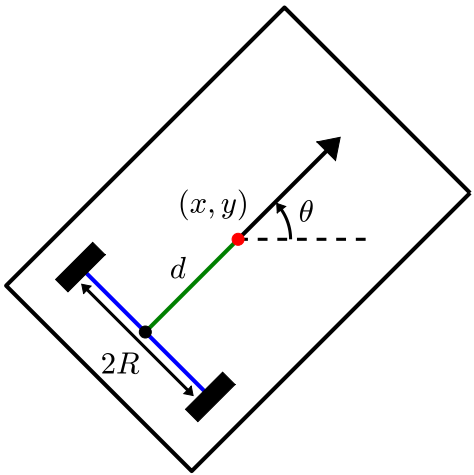}
\caption{A simple rectangular car.}
\label{fig:car}
\end{figure} 

As automated driving technology becomes more prevalent, it is ever more important to develop interpretable trajectory planning algorithms. In this manuscript, we address the problem of trajectory planning for simple self-driving cars using a method rooted in optimal control theory and dynamic programming. We consider the vehicle pictured in \cref{fig:car}. The configuration space for the car is $(x,y,\theta)$ where $(x,y)$ denotes the coordinate for the center of mass of the car, and $\theta \in [0,2\pi)$ denotes the angle of inclination from the horizontal. The rear axle has length $2R$ and the distance from the rear axle to the center of mass is $d$. Such cars are typically propelled using actuators which supply torque to either of the rear wheels \cite{Bertozzi}. The motion is subject to a nonholonomic constraint \begin{equation} 
\dot y \cos \theta - \dot x \sin \theta = d\dot \theta.
\end{equation} This ensures motion (approximately) tangential to the rear wheels; indeed, in the case that $d = 0$, the constraint reduces to $dy/dx = \tan \theta$. We also assume the car has a minimum turning radius, or equivalently, a maximum angular velocity so that $\lvert \dot \theta\rvert \le W$, for some $W >0$. LaValle \cite[Chap. 13]{LaValle} includes an extended discussion of models for this and similar vehicles.

Trajectory planning for simple cars goes back to Dubins \cite{Dubins} who considered that case that $d = R =0$ (so that the car is simplified to a point mass) and only allowed unidirectional (``foward'') movement. Reeds and Shepp \cite{ReedsShepp} considered forward and backward motion, and proved that in the absence of obstacles, the optimal trajectories are combinations of straight lines and arcs of circles of minimum radius, and that optimal trajectories have at most two kinks where the car changes from moving forward to backward or vice versa. Later effort was devoted toward adding obstacles \cite{Barraquand}, and developing an algorithm for near-optimal trajectories which are robust to perturbation \cite{AgarwalWang}. All of this work was carried out in a discrete and combinatorial fashion, breaking the paths into ``turning'' or ``straight'' segments and proving results regarding the possible combinations of these pieces. 

To the authors' knowledge, this problem was first analyzed in the context of optimal control theory by Boissonnat et al. \cite{boissonnat1,boissonnat2,boissonnat3} who gave shorter proofs and extensions of results of \cite{Dubins} and \cite{ReedsShepp}. Later, Takei and Tsai et al. \cite{TakeiTsai1,TakeiTsai2} used dynamic programming to derive a partial differential equation (PDE) which is solved by the optimal travel time function. Through all this work, the car was still simplified to a point mass. Later, the same approach was applied while considering the rectangular vehicle pictured in \cref{fig:car} \cite{ParkinsonCar}. PDE-based optimal path planning algorithms have also been developed for a number of applications besides simple self-driving cars, including underwater path planning in dynamic currents \cite{Lolla}, human navigation in a number of contexts, \cite{cartee2019time,Parkinson,Parkinson2} and recent models for environmental crime \cite{Arnold,cartee2020control,Chen}.

Other recent work has been devoted to machine learning and variatial approaches to the problem; for example, \cite{Shukla,Gao,Johnson}. Such approaches often rely on a hierarchical algorithms with global trajectory generation and local collision avoidance as in \cite{Rodriguez,Mao}. 

\subsection{Our Contribution}

We present a PDE based optimal path planning algorithm for simple self-driving vehicles. Our method is in the same spirit as \cite{TakeiTsai1,TakeiTsai2,ParkinsonCar}. We use dynamic programming to derive a Hamilton-Jacobi-Bellman (HJB) equation which is satisfied by the optimal travel time function. The optimal steering plan is generated using the solution to the HJB equation. To the authors' knowledge, in all previous HJB formulations of optimal trajectory planning for simple self-driving cars, the obstacles are stationary. We present a time-dependent formulation in which obstacles are allowed to move, which is a significant step in adding realism to this formulation. 

In general, the time depedent HJB equation has the form \begin{equation} \label{eq:generalHJB }u_t + H(\bx, \nabla u (\bx)) = r(x).\end{equation} Because the equation is nonlinear and solutions develop kinks \cite{Falcone}, some care is needed when solving HJB equations numerically; for example, they are not amenable to the simplest finite difference methods. Accordingly, we present an upwind finite difference scheme to solve our equation. 

The HJB formulation of minimal time path planning has a number of natural advantages. Because it is rooted in optimal control theory, there are some theoretical guarantees and there are no ``black box'' components, so the results are interpretable. It is also a very robust modeling framework, wherein one can easily account for a number of other realistic concerns such as energy minimization. Finally, it eschews the need for hierarchical algorithms, and after a single PDE solve, this formuation can resolve optimal paths from any initial configuration to the desired ending configuration.

\section{MATHEMATICAL FORMULATION}

Our algorithm is based on a control theoretic formulation. Generally, to analyze control problems using dynamic programming, one derives a Hamilton-Jacobi-Bellman equation which is satisfied by the optimal travel time function. Solving the equation provides the optimal travel time from any given starting configuration to a fixed ending configuration, and the derivatives of the travel time function determine the optimal steering plan. For general treatment of this approach (in both theory and practice), see \cite{Fleming,Bertsekas}.

\subsection{Equations of Motion \& Control Problem}

We consider a kinematic model of a self driving car which moves about a domain $\Omega \subset \R^2$ in the presence of moving obstacles. Fix a horizon time $T > 0$. At any time $t \in [0,T]$, the obstacles occupy a set $\Omega_{\text{obs}}(t) \subset \Omega$, so that the free space is given by $\Omega_{\text{free}}(t) = \Omega \setminus \Omega_{\text{obs}}(t)$. 

As described above, we track the current configuration of the car using variables $(\bx,\by,\bs) : [0,T] \to \Omega \times [0,2\pi)$, which obey the kinematic equations \begin{equation} \label{eq:kinematic}
\begin{split}
\dot \bx &= \bv \cos \bs - \bw W d \sin\bs, \\
\dot \by &= \bv\sin \bs + \bw W d \cos \bs, \\
\dot \bs &= \bw W. 
\end{split}
\end{equation}
Here $W > 0$ is a bound on the angular velocity of the vehicle which enforces bounded curvature of the trajectory. The control variables are $\bv(\cdot),\bw(\cdot) \in [-1,1]$, representing the tangential and angular velocity respectively. By taking velocities to be control variables, we are neglecting some of the ambient dynamics. In a more complete dynamic model, the control variables would be the torques applied to the rear wheels by the actuators. For a derivation of both this kinematic model and a dynamic model, see \cite{Moret}, and for generalizations of the kinematic model, see \cite{Triggs}.

For configurations $(x,y,\theta) \in \Omega \times [0,2\pi)$, define $C(x,y,\theta) \subset \Omega$ to be the space occupied by the car. In general, this could be any shape but for our purposes it will be a rectangle as pictured in \cref{fig:car}. Then given a desired ending configuration $(x_f,y_f,\theta_f)$, a trajectory $(\bx(t),\by(t),\bs(t))$ is referred to as \emph{admissable} if each of the following is true: \begin{itemize}
\item[(1)] it obeys \eqref{eq:kinematic} for $t \in [0,T]$,
\item[(2)] $C(\bx(t),\by(t),\bs(t)) \cap \Omega_{\text{obs}}(t) = \varnothing$ for all $t \in [0,T]$,
\item[(3)] $(\bx(T),\by(T),\bs(T)) = (x_f,y_f,\theta_f)$.
\end{itemize}
Here (2) signifies that the car does not collide with obstacles, and (3) signifies that the trajectory ends at the desired ending configuration. 

Given an initial configuration $(x,y,\theta)$, the goal is then to resolve the steering plan $\bv(t),\bw(t)$ that determines the minimal time required to traverse an admissable trajectory from $(x,y,\theta)$ to $(x_f,y_f,\theta_f)$. 

\subsection{The Dynamic Programming Approach}

We resolve the optimal steering plan using dynamic programming and a Hamilton-Jacobi-Bellman (HJB) equation. To analyze the problem in the dynamic programming framework, we first define the travel-time function. For a given configuration $(x,y,\theta) \in \Omega \times [0,2\pi)$ and time $t \in [0,T]$, we restrict ourselves to trajectories $(\bx(\cdot),\by(\cdot),\bs(\cdot))$ such that $(\bx(t),\by(t),\bs(t)) = (x,y,\theta)$. For such trajectories, if $\bv(\cdot),\bw(\cdot)$ is the corresponding steering plan, we define the first arrival time \begin{equation}\label{eq:arrivalTime} t^*_{\bv,\bw} = \inf\{s \, : \, (\bx(s),\by(s),\bs(s)) = (x_f,y_f,\theta_f)\}.\end{equation} The cost functional for the control problem is then \begin{equation}\label{eq:functional} \mathscr{T}(x,y,\theta,t,\bv(\cdot),\bw(\cdot)) = \left\{\begin{matrix}t^*_{\bv,\bw}, & \text{ if } t^*_{\bv,\bw} \le T, \\ +\infty, & \text{otherwise}. \end{matrix} \right.\end{equation} The optimal travel time function is then defined \begin{equation} \label{eq:valueFunc}u(x,y,\theta,t) = \inf_{\bv(\cdot),\bw(\cdot)} \mathscr{T}(x,y,\theta,t,\bv(\cdot),\bw(\cdot)). \end{equation} Intuitively, $u(x,y,\theta,t)$ is the minimal time required to steer the car to $(x_f,y_f,\theta_f)$, given that the car is at $(x,y,\theta)$ at time $t$. Note that if $(x,y,\theta)$ is far from $(x_f,y_f,\theta_f)$ and $t$ is close to $T$, there may be no way to steer the car to the ending configuration in the allotted time. If this is the case, then $u(x,y,\theta,t) = +\infty$. However, if there are any admissable trajectories $(\bx(\cdot),\by(\cdot),\bs(\cdot))$ such that $(\bx(t),\by(t),\bs(t)) = (x,y,\theta)$, then $u(x,y,\theta,t) \le T$. 

We want to derive a partial differental equation satisfied by the optimal travel time function. The dynamic programming principle \cite{Bellman} for this control problem is \begin{equation} \label{eq:DPP} \begin{split}u(&x,y,\theta,t) = \\ &\delta + \inf_{\bv(\cdot),\bw(\cdot)}\{u(\bx(t+\delta), \by(t+\delta), \bs(t+\delta),t+\delta)\}\end{split} \end{equation} where $(\bx(t),\by(t),\bs(t)) = (x,y,\theta)$ and the infimum is taken with respect to the values $\bv(s),\bw(s)$ for $s \in (t,t+\delta)$. 

Supposing that $u(x,y,\theta,t)$ is smooth, we can divide by $\delta$ and take the limit as $\delta \to 0$ to arrive at \begin{equation} \label{eq:preHJB} \inf_{v,\w} \big\{ u_t + \dot \bx u_x + \dot \by u_y + \dot \bs u_{\theta} \big\} = -1, \end{equation} whereupon inserting \eqref{eq:kinematic} yields \begin{equation}
u_t + \inf_{v,\w} \left\{ \begin{split} (&u_x \cos\theta + u_y \sin \theta)v\,\,+ \\  &W(-du_x\sin\theta + du_y\cos\theta+ u_\theta)\w \end{split}\right\} = -1.
\end{equation} Notice the minimization is linear in the variables $v,\w \in [-1,1]$, and thus the minimizing values can be resolved explicitly. We see that \begin{equation} \label{eq:controls}
\begin{split} v &= - \text{sign}(u_x \cos\theta + u_y \sin \theta), \\ \w &= -\text{sign}(-d\sin\theta u_x + d\cos\theta u_y + u_\theta), \end{split}
\end{equation} where $u(x,y,\theta,t)$ solves the HJB equation \begin{equation} 
\label{eq:HJB} 
\begin{split} 
u_t& - \abs{u_x \cos \theta + u_y \sin \theta}\\ &\,\,\,\,\,\,- W\abs{-d u_x \sin\theta + du_y \cos \theta + u_\theta}
\end{split} = -1.
\end{equation} This derivation is only valid when $u(x,y,\theta,t)$ is smooth, which is not expected to be the case. However, under very general conditions, the travel time function is the unique viscosity solution of \eqref{eq:HJB}  \cite{CrandallLions}. For a fully rigorous derivation of the Hamilton-Jacobi-Bellman equation, see \cite{Bardi1997}.

There are a few natural conditions appended to \eqref{eq:HJB}. At the terminal time $T$, the cost functional \eqref{eq:functional} assigns a value of either 0 or $+\infty$, depending on whether car is at the ending configuration or not. Thus, we have the terminal condition \begin{equation}\label{eq:terminal}
u(x,y,\theta,T) = \left\{\begin{matrix} 
0, & (x,y,\theta) = (x_f,y_f,\theta_f), \\ +\infty, & \text{otherwise},
\end{matrix} \right.
\end{equation} and we want to resolve $u(x,y,\theta,t)$ for \emph{preceding} times $t \in [0,T)$. So the equation runs ``backwards'' in time.

Likewise, if the trajectory has already arrived at the ending configuration, the remaining travel time is 0, so we have the boundary condition \begin{equation}\label{eq:boundary}
u(x_f,y_f,\theta_f,t) = 0, \,\,\,\,\,\, t \in [0,T].
\end{equation} Lastly, to ensure the car does not collide with obstacles, we assign $u(x,y,\theta,t) = +\infty$ for any $(x,y,\theta,t)$ such that $C(x,y,\theta) \cap \Omega_{\text{obs}}(t) \neq \varnothing$. 

By \eqref{eq:controls}, the only possible values of the control variables are $v,\w\in\{-1,0,1\}$, resulting in a bang-bang controller which has a ``no bang'' option.  This makes intuitive sense because there is never incentive to drive or turn slower than the maximum possible speed, unless one needs to wait for an obstacle to move out of the way (whereupon $v=0$) or one needs to drive in a straight line (whereupon $\w = 0$). When no obstacles are present, one can eliminate the $v = 0$ option and the path will consist of straight lines and arcs of circles of minimum radius, which agrees with early analysis of the problem \cite{Dubins,ReedsShepp}.

As a final note, this derivation is very similar to that in \cite{ParkinsonCar}. However, when the obstacles are stationary, as in \cite{TakeiTsai1,TakeiTsai2,ParkinsonCar}, the optimal travel time function does not depend on $t$, since the optimal trajectory depends only upon the current configuration, not upon the time $t$ when the car occupies that configuration. In that case, one can eliminate the time horizon $T$, and opt instead for a stationary HJB equation. One can than visualize solving the stationary HJB equation by evolving a front outward from the final configuration, and recording the time as the front passes through other configurations, terminating when each point in the domain $\Omega \times [0,2\pi)$ has been assigned a value. This is the philosophy behind level-set inspired optimal path planning \cite{Parkinson,Parkinson2}, and numerical implementations like fast sweeping \cite{tsai2003fast,kao2004lax,RotatingGrid} and fast marching methods \cite{Tsitsiklis,Sethian1}. In theory, something similar is possible here. If one does not care to enforce a finite time horizon, then making the substitution $\tau = T-t$ and taking $T \to \infty$ will do away with it. However, in practice, we will want to discretize the HJB equation in order to solve computationally, which will require choosing a fixed time horizon. Thus we cannot do away with $T$, but to minimize its effect, we set it large enough that the travel time $u(x,y,\theta,0)$ is finite for all $(x,y,\theta) \in \Omega \times [0,2\pi)$. In this manner, any initial configuration (not overlapping the obstacles) will have admissable paths which reach $(x_f,y_f,\theta_f)$ within time $T$.

\section{NUMERICAL METHODS}

In this section, we design a numerical scheme to approximate \eqref{eq:HJB}. Since Hamilton-Jacobi equation admit non-smooth solutions which cannot be approximated by simple finite difference schemes, effort has been expended to develop schemes which resolve the viscosity solution. For a survey of numerical methods for Hamilton-Jacobi equations, see \cite{Falcone,osher2003level}.

\subsection{An Upwind, Monotone Scheme for \eqref{eq:HJB}}

For simplicity, we confine ourselves to a rectangular spatial domain $\Omega = [x_{\text{min}}, x_{\text{max}}] \times [y_{\text{min}},y_{\text{max}}]$. Choosing $I,J,K,N \in \mathbb N$, let $(x_i)_{i=0}^I, (y_j)_{j=0}^{J}, (\theta_k)_{k=0}^K, (t_n)_{n=0}^N$ be uniform discretizations of their respective domains with grid parameters $\Delta x, \Delta y, \Delta \theta, \Delta t$, and let $u_{ijk}^n$ be our approximation to $u(x_i,y_j,\theta_k,t_n)$. For each $v,\w \in \{-1,0,1\}$, define \begin{equation} \label{eq:coeffs} 
\begin{split}
A_k(v,\w) &= v\cos \theta_k - \w Wd\sin\theta_k, \\
a_k(v,\w) &= \text{sign}(v\cos \theta_k - \w Wd\sin\theta_k), \\
B_k(v,\w) &= v\sin \theta_k + \w W d \cos \theta_k,\\
b_k(v,\w) &= \text{sign}(v\sin \theta_k + \w W d \cos \theta_k).
\end{split}
\end{equation} Then \eqref{eq:preHJB} can be rewritten \begin{equation}\label{eq:discEq}
u_t + \min_{v,\w}\{A_k(v,\w)u_x +B_k(v,\w)u_y + \w W u_\theta\} = -1.
\end{equation} 
Recall, the terminal values $u^N_{ijk}$ are supplied here, and we need to integrate this equation backwards in time. Thus at time step $t_n$, we need to resolve $u^n_{ijk}$ given known values $u^{n+1}_{ijk}$. This suggests backward Euler time integration \begin{equation} \label{eq:timeDisc} (u_t)^n_{ijk} = \frac{u^{n+1}_{ijk} - u^{n}_{ijk}}{\Delta t}. \end{equation} The upwind approximations to the other derivatives in \eqref{eq:discEq} using $u^{n+1}_{ijk}$ are given by 
\begin{equation} \label{eq:upwindApprox}
\begin{split}
(A_k(v,\w)u_x)^{n+1}_{ijk} &= \abs{A_{k}(v,\w)} \left(\frac{u^{n+1}_{i+a_k(v,\w),j,k} - u^{n+1}_{ijk}}{\Delta x}\right), \\
(B_k(v,\w)u_y)^{n+1}_{ijk} &= \abs{B_{k}(v,\w)}\left( \frac{u^{n+1}_{i,j+b_k(v,\w),k} - u^{n+1}_{ijk}}{\Delta y}\right), \\
(\w W u_\theta)^{n+1}_{ijk} &= \abs \w W \left(\frac{u^{n+1}_{i,j,k+\text{sign}(\w)} - u^{n+1}_{ijk}}{\Delta\theta}\right).
\end{split}
\end{equation}
We insert these approximations in \eqref{eq:discEq} to arrive at \begin{equation} \label{eq:update} \begin{split} u^n_{ijk} = u^{n+1}_{ijk} + &\Delta t \left(1 + \min_{v,w} \{(A_k(v,\w)u_x)^{n+1}_{ijk} \right.\\&+ \left.(B_k(v,\w)u_y)^{n+1}_{ijk} + (\w W u_\theta)^{n+1}_{ijk} \} \right). \end{split} \end{equation} Since there are only finitely many pairs $(v,\w)$, we can compute the right hand side for each pair and explicity choose the pair which suggests the minimum possible value. Using this formula and stepping through $n=N-1, N-2,\ldots, 1,0$, we arrive at our approximation of the travel time function. 

To initialize, we set $u_{ijk}^n = +\infty$ (or some very large number) for all $i,j,k,n$ except at the node $(i_f,j_f,k_f)$ respresenting the configuration nearest to $(x_f,y_f,\theta_f)$ where we set $u^n_{i_f,j_f,k_f} = 0$ for all $n$. We then only update the node $u^n_{ijk}$ if the value suggested by \eqref{eq:update} is smaller than the value already stored at $u^n_{ijk}$. This ensures that the scheme is monotone so long as the CFL condition \begin{equation} \label{eq:CFL} 
\Delta t \left(\frac{1+Wd}{\Delta x} + \frac{1+Wd}{\Delta y} + \frac{W}{\Delta \theta} \right) \le 1
\end{equation} is satisfied \cite{Falcone,osher2003level}. In this case, since the scheme is also consistent, the approximation converges to the viscosity solution of \eqref{eq:HJB} as $\Delta x, \Delta y, \Delta \theta, \Delta t \to 0$.

We include a few implementation notes. First, to account for obstacles, at each time step $n$, we first need to find the illegal nodes (i.e. those which correspond to configurations wherein the car collides with an obstacle). At these nodes $(i^*,j^*,k^*)$, we do not use \eqref{eq:update}, but rather set $u^n_{i^*,j^*,k^*} = +\infty$. In previous work, this could be done in pre-processing since the obstacles were stationary and illegal configurations only needed to be resolved once. In this work, since the obstacles move, this must be repeated at every time step. Second, we use \eqref{eq:update} for $i = 2,\ldots,I-1, j = 2,\ldots,J-1$. The values $u_{ijk}^n$ at nodes corresponding to the spatial boundary are never updated, but should be given the value $+\infty$. This will ensure that the car never leaves the domain. Because we enforce the correct causality, the boundary nodes have no effect on interior nodes. Third, one needs to enforce periodic boundary conditions in $\theta$ by identifyting the nodes at $k=0$ and $k=K$. Lastly, above it is stated that $v,\w \in \{-1,0,1\}$. However, because it is impossible to turn a car without moving backward or forward, one should eliminate the cases $(v,\w) = (0,\pm1)$. So there are seven possible pairs of $(v,\w)$ to consider in total; in short, $(v,\w) = (\pm1, \pm1), (\pm1, 0), (0,0)$.

\subsection{Generating Optimal Trajectories}

There are a few different manners in which one can obtain optimal control values and generate optimal trajectories. It is possible to resolve control values $v^n_{ijk},\w^n_{ijk}$ while evaluating \eqref{eq:update}. One can define them to be the pair that achieves the minimum in \eqref{eq:update} at any node $(i,j,k,n).$ Alternatively, after resolving $u^{n}_{ijk}$, one can interpolate to off-grid values and use \eqref{eq:controls} to resolve the optimal steering plan at any point $(x,y,\theta,t)$. In either case, after choosing an initial point, one can insert the optimal control values into \eqref{eq:kinematic} and integrate the equations of motion until the trajectory reaches $(x_f,y_f,\theta_f)$. This is the approach taken by \cite{ParkinsonCar}.

In a different approach, we opt for a semi-Lagrangian path-planner as in \cite{TakeiTsai2,cartee2019time}. Specifically, we first interpolate $u^n_{ijk}$ to off grid values, so we have an approximate travel time function $u(x,y,\theta,t)$. Then, choosing an initial point $(\bx_0,\by_0,\bs_0)$ and a time step $\delta > 0$, and rewriting \eqref{eq:kinematic} as $(\dot \bx, \dot \by,\dot \bs) = F(\bx,\by,\bs,\bv,\bw)$, we set \begin{equation} \label{eq:SLscheme}
\begin{split}
(v^*,\w^*) = \argmin_{v,\w} u( (\bx_\ell,\by_\ell,\bs_\ell) + \delta F(\bx_\ell,\by_\ell,\bs_\ell,v,\w),\ell\delta),\\
(\bx_{\ell+1},\by_{\ell+1},\bs_{\ell+1}) = (\bx_\ell,\by_\ell,\bs_\ell) + \delta F(\bx_\ell,\by_\ell,\bs_\ell,v^*,\w^*),
\end{split}
\end{equation} for $\ell = 0,1,2,\ldots$, halting when $(\bx_\ell,\by_\ell,\bs_\ell)$ is within some tolerance of $(x_f,y_f,\theta_f)$.

\section{RESULTS \& EXAMPLES} \label{sec:examp}

We present results of our algorithm in three examples. In all cases, we use the spatial domain $\Omega = [-1,1] \times [-1,1]$.  We take the car to be a rectangles as pictured in \cref{fig:car} with $d = 0.07$ and $R = 0.04$ and we take the maximum angular velocity to be $W = 4$. These are dimensionaless variables used for testing purposes. In each of the following pictures the final configuration $(x_f,y_f,\theta_f)$ will be marked with a red star and the initial configurations of the various cars will be marked with green stars. We use a $101 \times 101 \times 101$ discretization of $\Omega \times [0,2\pi)$ and then choose $\Delta t$ according to the CFL condition \eqref{eq:CFL}. We choose the time horizon $T= 10$. As mentioned before, this simply needs to be chosen so that there are admissable paths from every point on the domain to the final configuration which take time less than $T$ to traverse. In some of the examples it could likely be smaller, but $T = 10$ was sufficiently large for all of them.

\begin{figure*}
\centering
\includegraphics[width=0.23\textwidth]{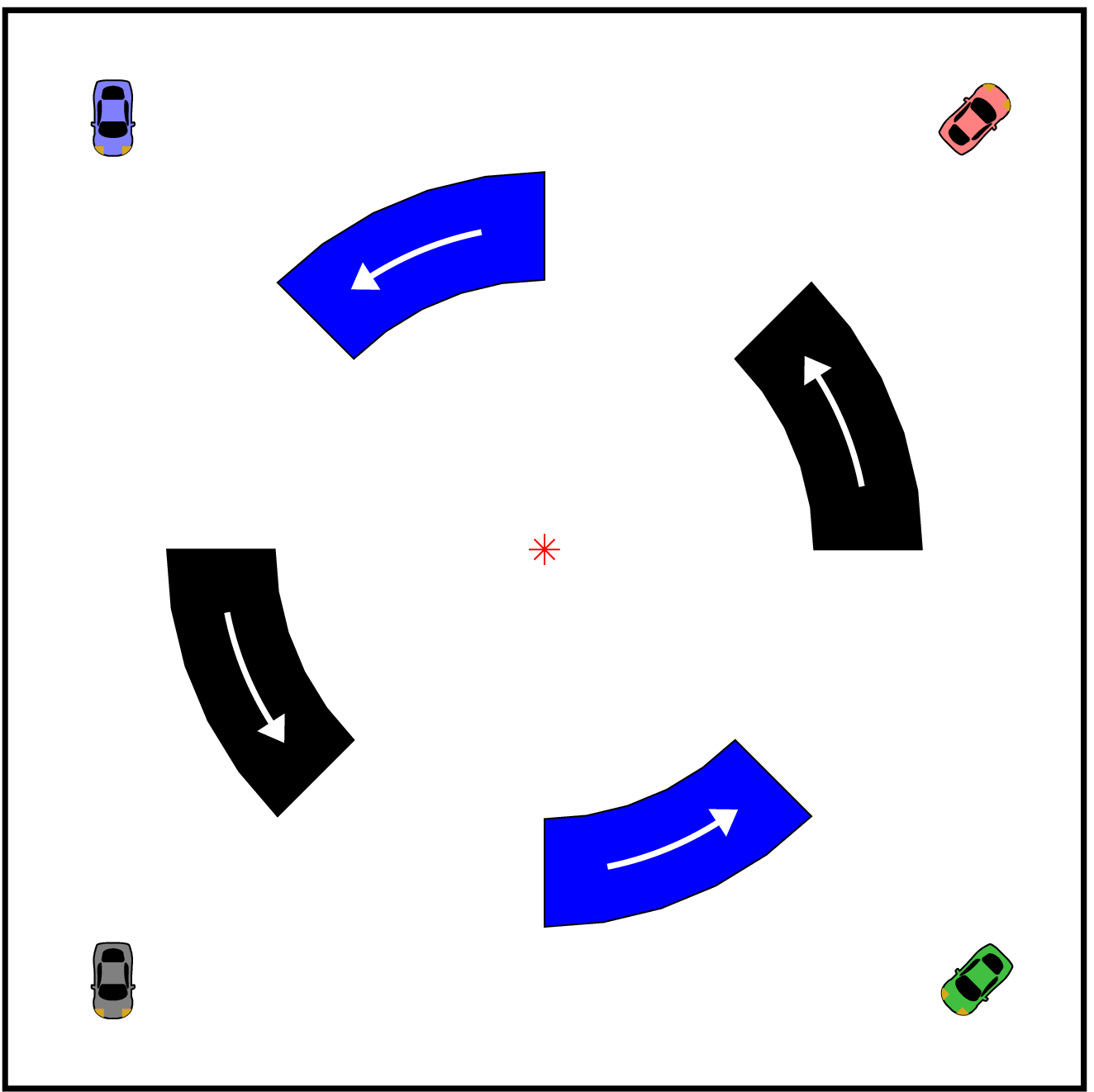} \,
\includegraphics[width=0.23\textwidth]{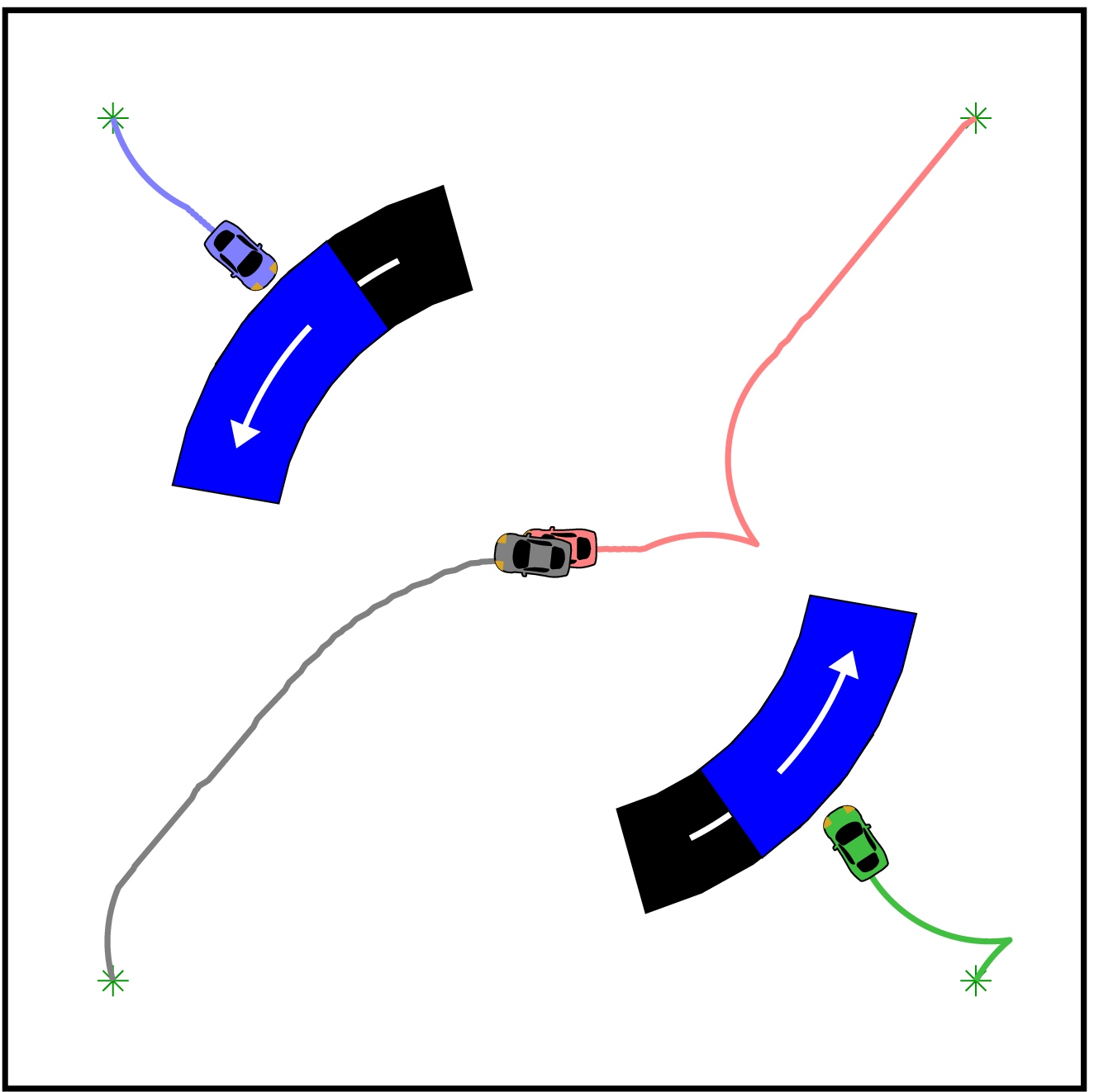} \,
\includegraphics[width=0.23\textwidth]{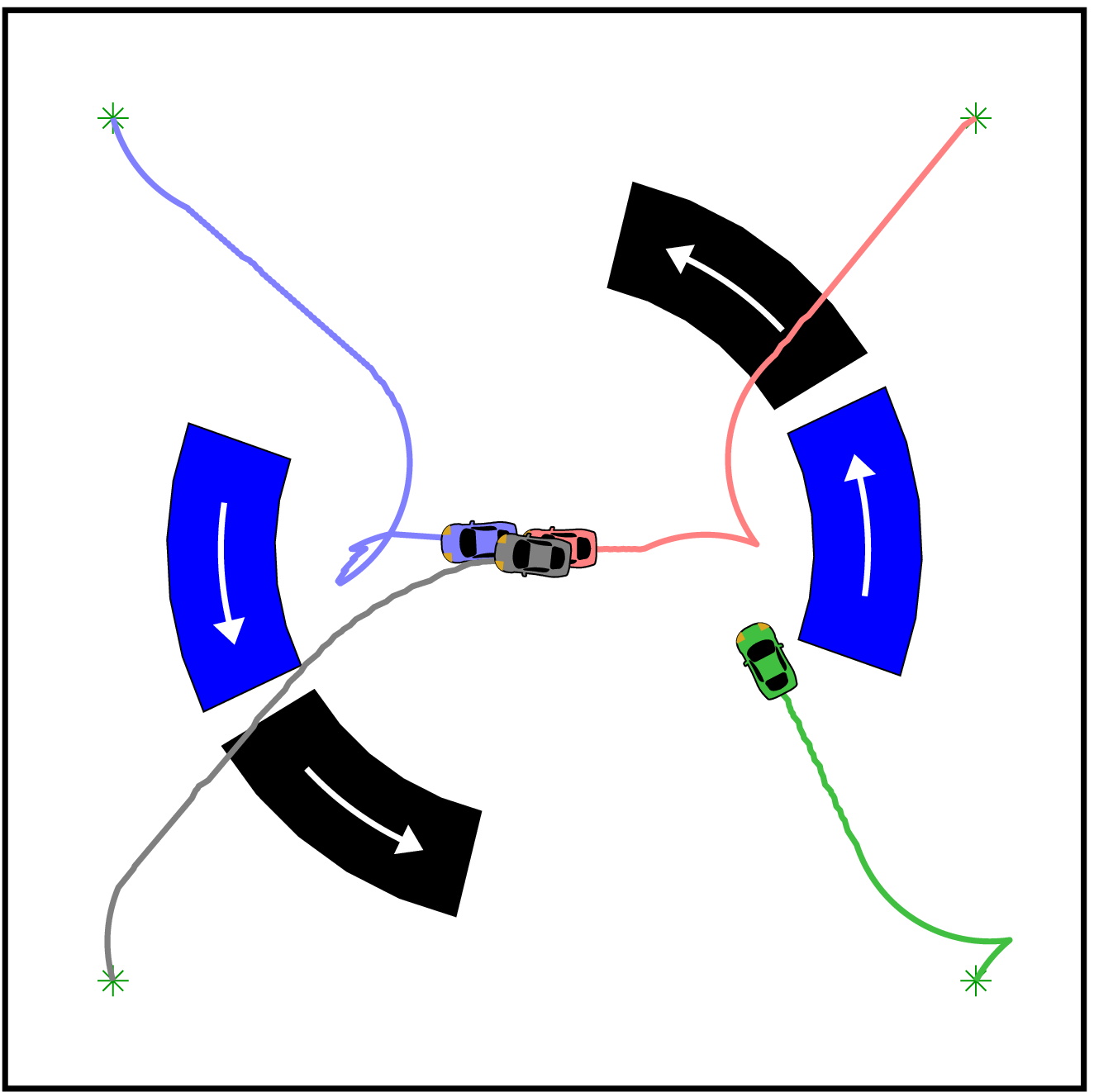} \,
\includegraphics[width=0.23\textwidth]{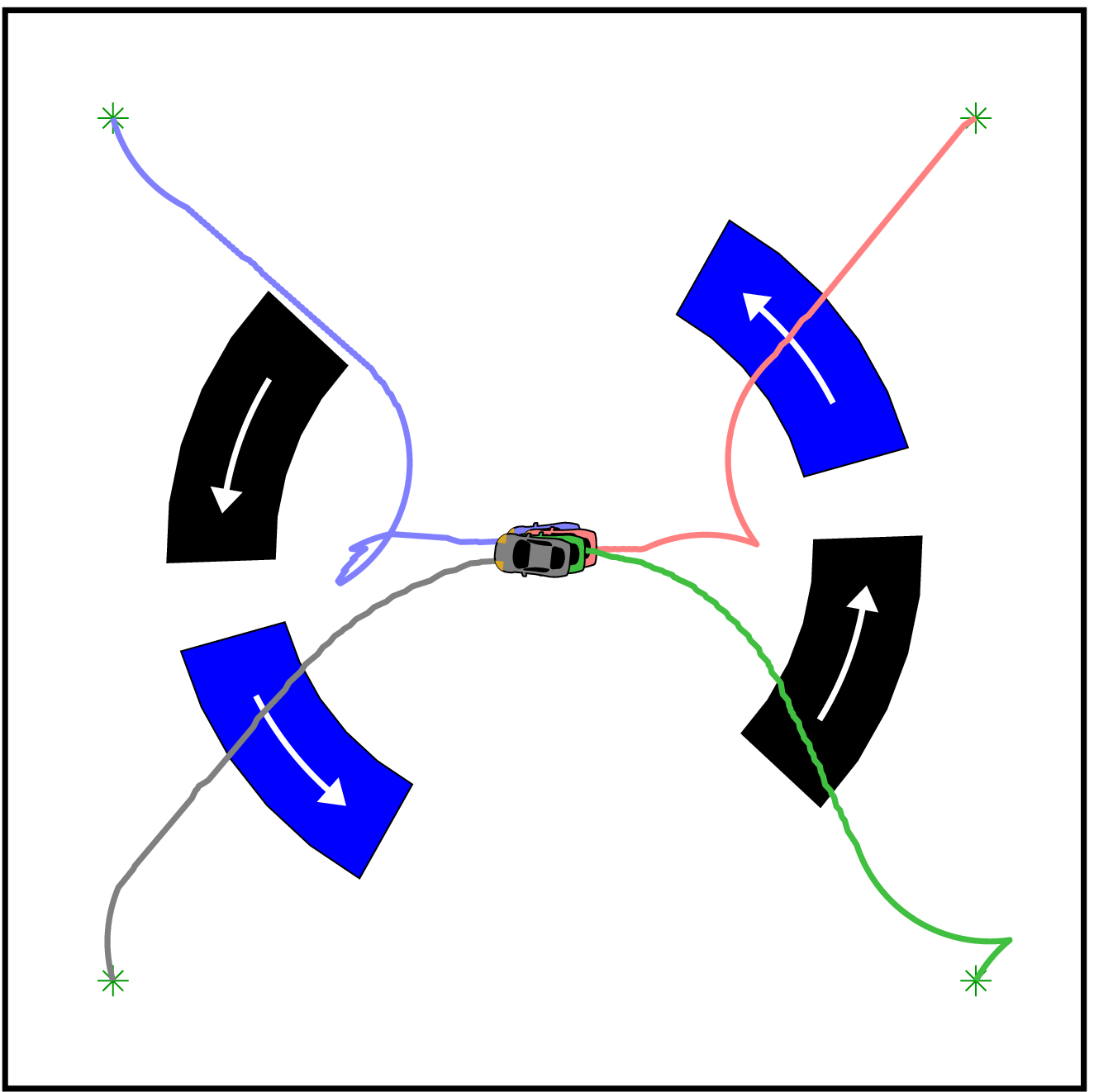}
\caption{In the first example, we have several cars starting in the corners and ending in the center facing due west. The obstacles (black and blue sectors) rotate counterclockwise with black obstacles rotating at three times the speed of the blue ones. }
\label{fig:examp1}
\end{figure*}

In the first example, the final configuration is $(x_f, y_f, \theta_f) = (0,0,\pi)$ meaning the cars will end at the center of the domain, facing due west. In this case, the obstacles $\Omega_{\text{obs}}(t)$ are 4 annular sectors which will rotate about the origin in the counterclockwise direction. These are represented in black and blue in \cref{fig:examp1}. The black obstacles rotate with 3 times the speed of the blue obstacles. Notice that in the second panel, the green and blue cars (respectively bottom right and top left corners) need to stop and wait to let the obstacles pass before completing their path. The grey and pink cars are essentially unoccluded and can travel directly to the destination. We note that these paths are generated individually, and simply plotted over each other. There is no competition between the cars.

\begin{figure*}
\centering
\includegraphics[width=0.23\textwidth]{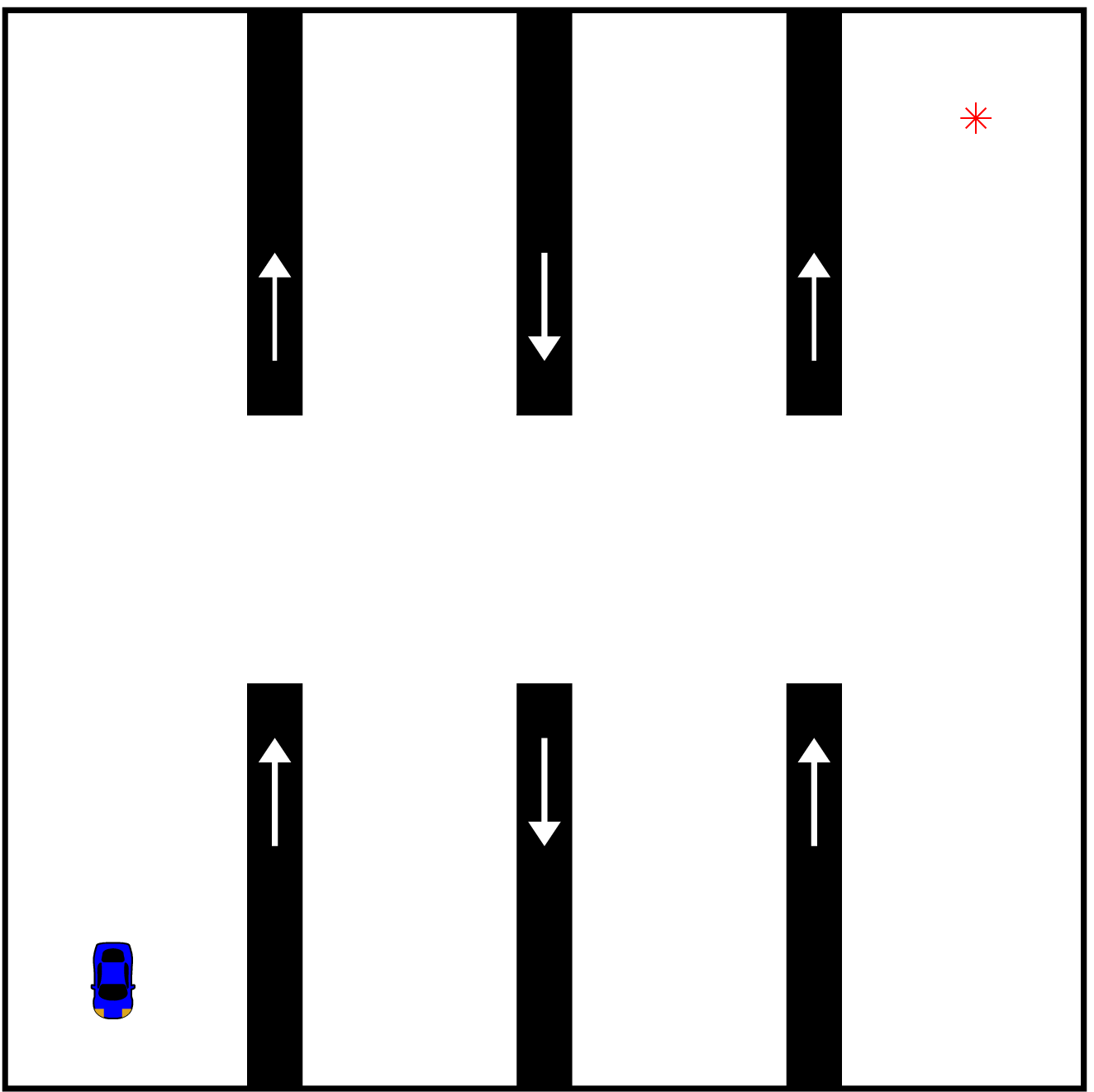} \,
\includegraphics[width=0.23\textwidth]{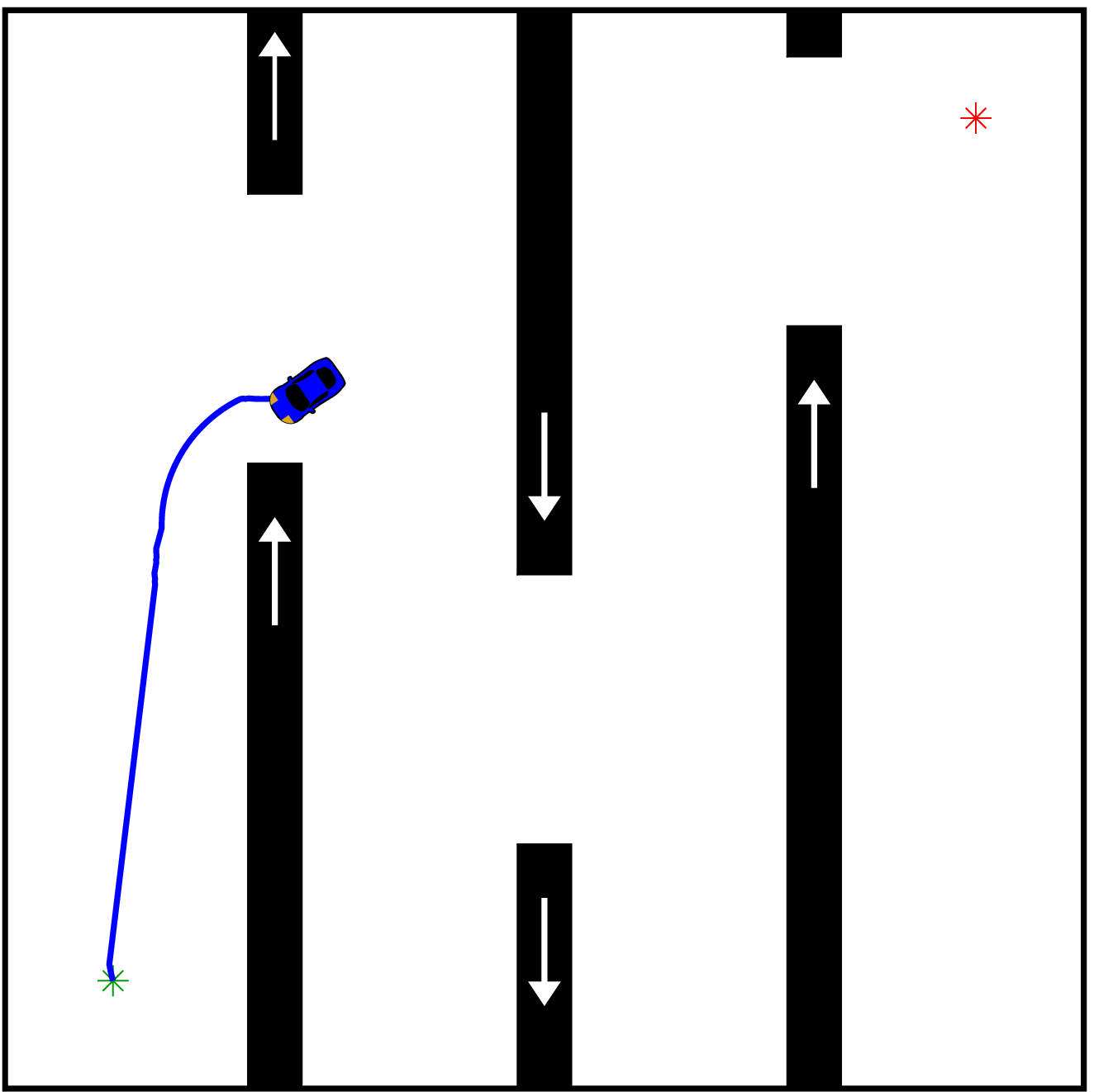} \,
\includegraphics[width=0.23\textwidth]{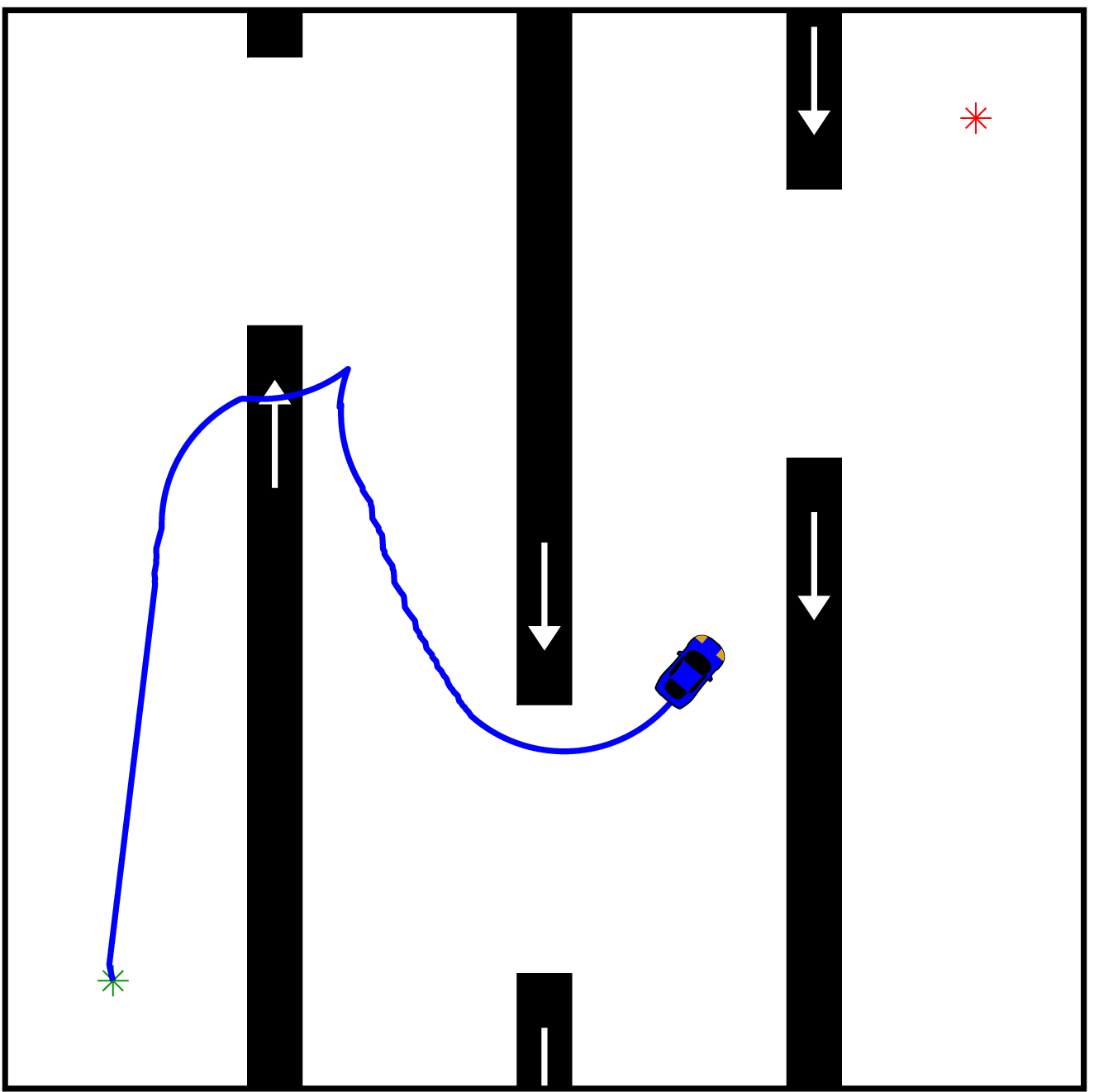} \,
\includegraphics[width=0.23\textwidth]{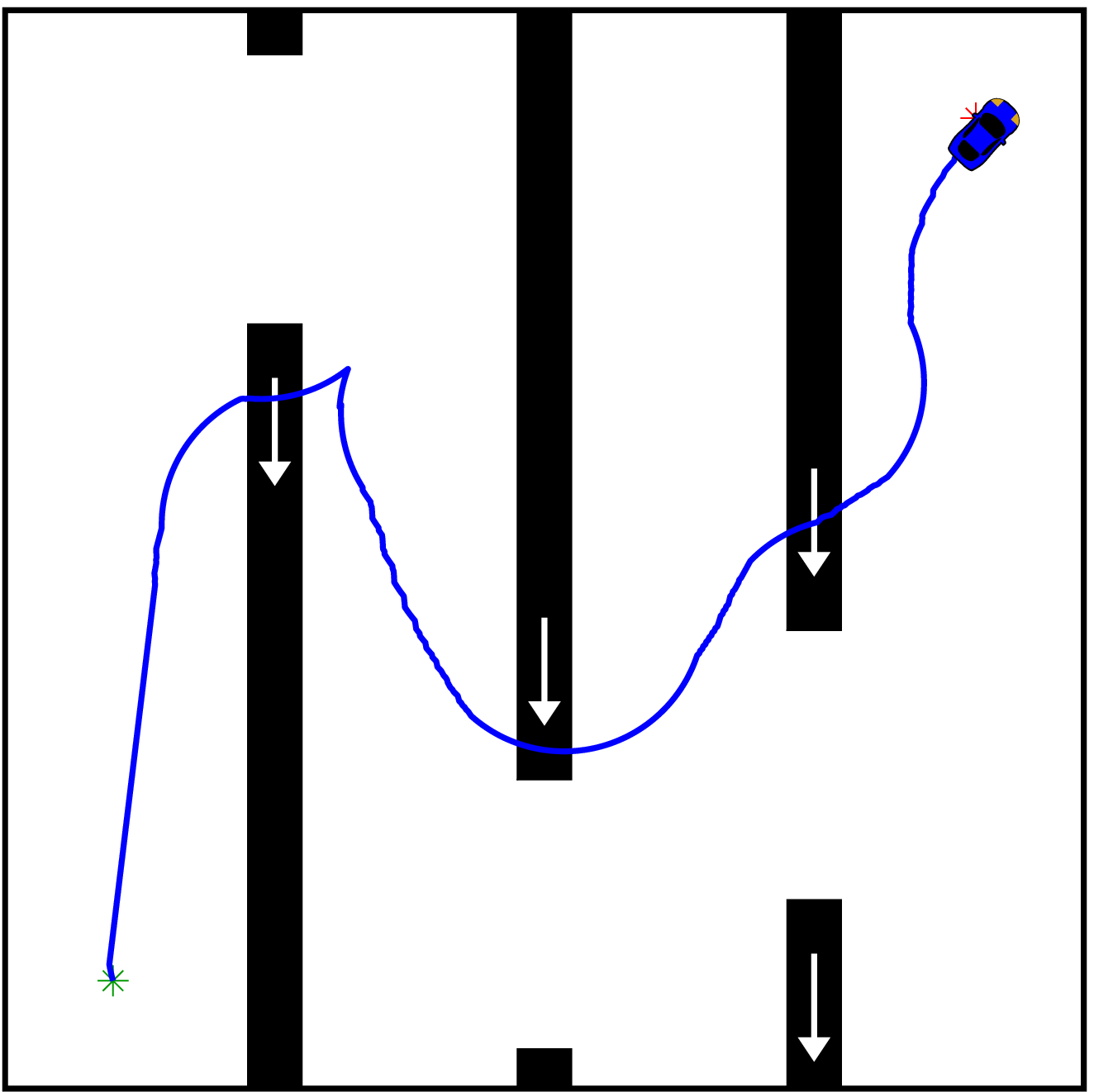}
\caption{In the second example, one care attempts to navigate through moving doorways. The black obstacles oscillate up and down as indicated by the arrows in each panel.}
\label{fig:doors}
\end{figure*}

In the second example, the final configuration is $(x_f,y_f,\theta_f) = (0.8, 0.8, \pi/4)$ so that the car needs to end near the top right corner of the domain facing northeast. The car begins in the bottom left corner of the domain as seen in \cref{fig:doors}, and must navigate through three moving doorways. The black bars represent the obstacles and they oscillate as indicated by the arrows. The car is able to navigate through the domain without stopping to wait for the doors. 

\begin{figure}
\centering
\includegraphics[width=0.1\textwidth]{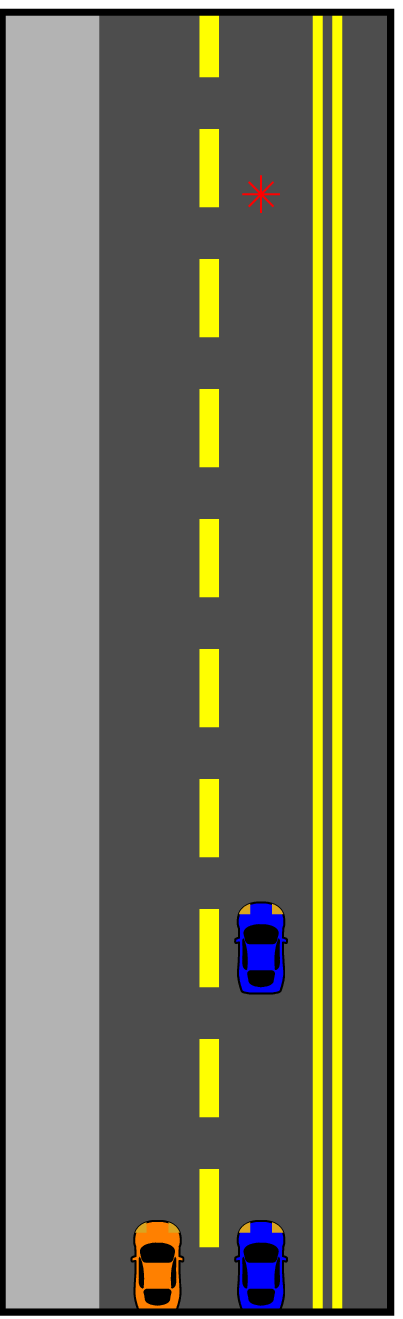} \,
\includegraphics[width=0.1\textwidth]{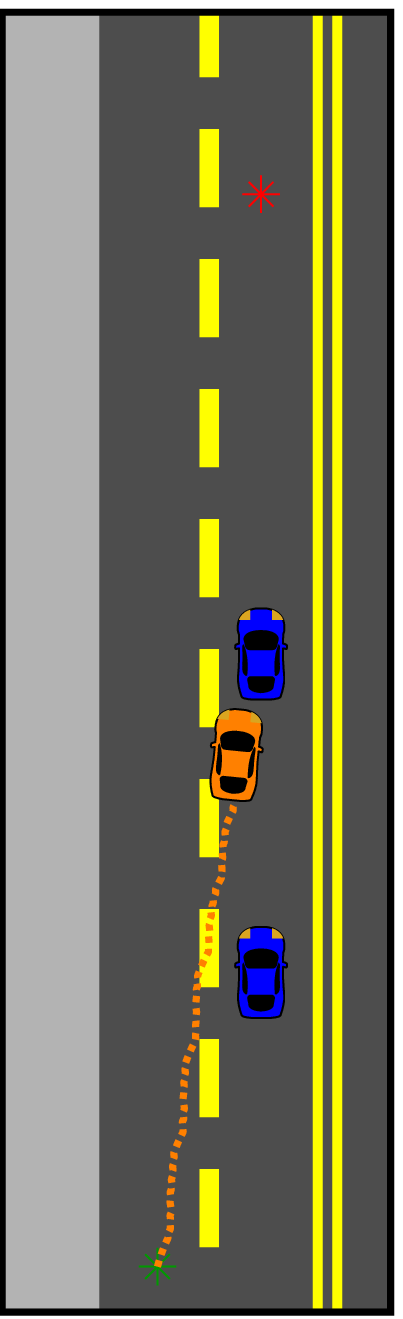} \,
\includegraphics[width=0.1\textwidth]{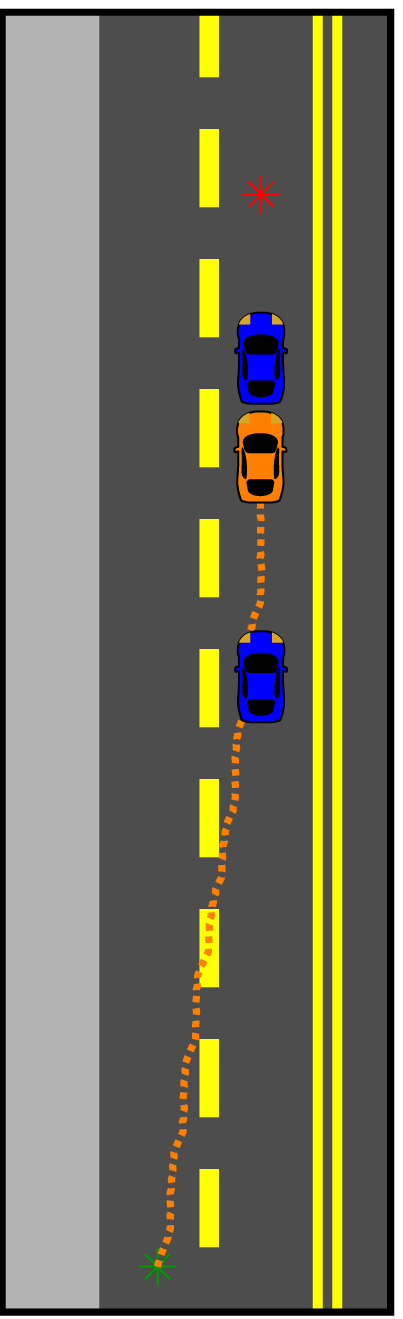} \,
\includegraphics[width=0.1\textwidth]{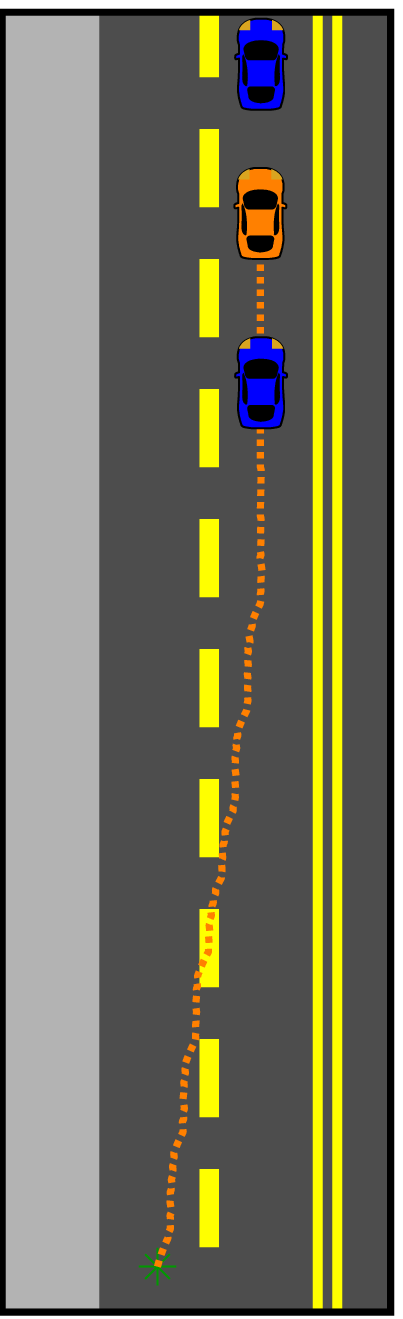}
\caption{A car (orange) changing lanes between two other cars (blue). Here the blue cars are the obstacles.}
\label{fig:changingLanes}
\end{figure}

In the third example, we consider the more realistic scenario of a car changing lanes in between two other cars as seen in \cref{fig:changingLanes}. In this case the two blue cars are treated as obstacles and the orange car must slide in between them.

\section{CONCLUSION \& DISCUSSION}

We present a Hamilton-Jacobi-Bellman formulation for time-optimal paths of simple vehicles in the presence of moving obstacles. This is distinguished from previous similar formulations which could only handle stationary obstacles. 

There are many ways in which this work could be extended. Some simple improvements would be to account for other realistic concerns such as energy minimization or intrumentation noise, which can both be added to the model in a straightforward manner, though they may complicate the numerical methods. 

Perhaps the biggest drawback of this method is that it is currently too computationally intensive for real-time applications. The simulations for each of the examples in \cref{sec:examp} required several minutes of CPU time (on the authors' home computers). However, one may be able to apply recent methods for high-dimensional Hamilton-Jacobi equations \cite{Darbon,Lin}. These methods are based on Hopf-Lax type formulas and trade finite differences for optimization problems. It may be difficult to account for crucial boundary conditions in our model when using such schemes, so some care would be required. However, if they could be applied to this problem, it would also provide an opportunity to extend the model to higher dimensions where finite difference methods are infeasible.

\section*{ACKNOWLEDGMENT}

The authors were supported in part by NSF DMS-1937229 through the Data Driven Discovery Research Training Group at the University of Arizona.

\bibliographystyle{ieeetr}
\bibliography{bibliography}

\end{document}